\documentclass[12pt]{amsart}
\usepackage{amssymb}
\newcommand{\nc}{\newcommand}
\nc{\sg}{\sigma} \nc{\al}{\alpha} \nc{\dt}{\delta} \nc{\bt}{\beta}
\nc{\tht}{\theta} \nc{\gm}{\gamma} \nc{\Gm}{\Gamma}
\nc{\vf}{\varphi} \nc{\ol}{\overline} \nc{\iy}{\infty}
\nc{\Tht}{\Theta} \nc{\pa}{\partial} \nc{\sbs}{\subset}
\nc{\eq}{\equiv} \nc{\bb}{1\hskip-0.9mm{\text I}}

\nc{\bC}{\mathbb{C}} \nc{\bT}{\mathbb{T}} \nc{\bo}{\mathbb{O}}

\newtheorem{theorem}{\bf \sc Theorem}

\nc{\const}{\operatorname{const}}

\begin{document}
\begin{center}
    {\bf On the Estimation of the Convergence Rate in the Janashia-Lagvilava Spectral Factorization Algorithm}

\vskip+0.5cm
   {Lasha Ephremidze and Nobuhiko Fujii }

\end{center}

\vskip+0.7cm

{\small {\bf Abstract.} In the present paper, we estimate the
convergence rate  in the Janashia-Lagvilava spectral factorization
algorithm (see Studia Mathematica, 137, 1999, 93-100) under the
restriction on a spectral density matrix that its inverse is
integrable.

\vskip+0.7cm

}

\footnotetext{ {}{\em 2000 Mathematics Subject Classification}\,:
47A68.

}

\vskip+0.7cm

\subsection{Introduction} Wiener's spectral factorization theorem
 [7], [8] (see also [2], [1]) asserts that a positive definite
matrix-valued function
\begin{equation}
    S(z)=(f_{ij}(z))_{i,j=1}^r,
\end{equation}
$f_{ij}\in L_1(\bT)$, on the unit circle of the complex plane with
integrable logarithm of the determinant
\begin{equation}
    \log\,\det\,(S(z))\in L_1(\bT),
\end{equation}
admits the representation
\begin{equation}
    S(z)=\chi^+(z)(\chi^+(z))^*,
\end{equation}
where
\begin{equation}
    \chi^+(z)= \sum_{k=0}^\iy\gamma_kz^k,
\end{equation}
$|z|<1$, $\gamma_k$ are $r\!\times\!r$ matrix coefficients, is an
analytic matrix-function with entries from the Hardy space $H_2$,
$\chi^+(z)\in H_2$, and the determinant of which is an outer
function.

The relation (3) is assumed to hold a.e. on the unit circle $\bT$
and $(\chi^+)^*=(\ol{\chi^+})^T$ is the adjoint of $\chi^+$.

The condition (2) is also necessary for the existence of spectral
factorization and with a suitable constraint on $\chi^+(0)$ we can
ensure that the spectral factor $\chi^+$ is unique.

It is widely known that the solution of several applied problems in
Control Engineering and Communications  require spectral
factorization.  Due to its importance, it is not surprising that,
starting from Wiener's efforts [8], a variety of methods has been
developed for approximate numerical calculation of coefficients
$\gamma_k$ in (4) for a given matrix-function (1) (see the survey
paper [6]). Most of these algorithms are applicable only under the
additional restriction on $S$, e.g., to be scalar or rational. There
exists an explicit formula for the spectral factor in the scalar
case,
\begin{equation}
\chi^+(z)=\exp\left(\frac 1{4\pi}
\int\nolimits_0^{2\pi}\frac{e^{it}+z}{e^{it}-z}\log
S(e^{it})\,dt\right),
\end{equation}
which is a core of Kolmogorov's exp-log method (see [4], p. 211).
There is no analog of formula (5) in the matrix case, since in
general $e^{A+B}\not=e^Ae^B$ for non-commutative matrices $A$ and
$B$. This is the main reason that the algorithms  in the matrix case
are significantly more difficult.

An absolutely new approach to the matrix spectral factorization
problem was proposed by G. Janashia and E. Lagvilava in [3] without
imposing on matrix-function (1) any additional restriction apart
from the necessary and sufficient condition (2). This is the first
time that  the theory of Hardy spaces is used for solution of the
problem which naturally turned out to be very effective since the
problem itself is posed in this branch of mathematics.

A sequence of matrix-functions $S_{n}$ is constructed in
Janashia-Lagvilava algorithm which approximates $S$ in $L_1$ norm,
$$
    \big\|S_{n}(z)-S(z)\big\|_{L_1}\to 0.
$$
Then the explicit spectral factorization of $S_{n}$ is performed
\begin{equation}
S_{n}(z)=\chi_{n}^+(z)(\chi_{n}^+(z))^*,
\end{equation}
 and it is proved that
\begin{equation}
    \big\|\chi_{n}^+(z)-\chi^+(z)\big\|_{H_2}\to 0.
\end{equation}

In the present paper, we intend to obtain some qualitative
estimation of the closeness of $\chi_{n}^+$ to $\chi^+$ which might
be useful for practical computations of an approximate spectral
factor within a given accuracy. We achieve this goal under the
certain restriction on (1), namely,
\begin{equation}
S^{-1}(z)\in L_1(\bT),
\end{equation}
and estimate
\begin{equation}
    \big\|\chi_{n}^+(z)-\chi^+(z)\big\|_{H_1}=\big\|\chi_{n}^+(z)-\chi^+(z)\big\|_{L_1}
\end{equation}
from above (see Theorem 1 below).
 For example, a trigonometric polynomial matrix-function
$S(z)=\sum_{k=-N}^N\sg_kz^k$, $\sg_k$ are $r\!\times\!r$ matrix
coefficients, without zeroes of the determinant on $\bT$, $\det
S(z)>0$, $|z|=1$, satisfies this condition. It is well-known that
the spectral factor $ \chi^+(z)= \sum_{k=0}^N\gamma_kz^k$ is
analytic polynomial of the same order $N$ in this case, and the
estimation of (9) can be used to obtain the accuracy of
approximately computed matrix coefficients $\gamma_k$,
$k=0,1,\ldots,N$.

In the present paper we consider only two dimensional matrices,
although, as Prof. Lagvilava informed us, the  method we propose can
be extended to higher dimensional matrices as well.

It should be mentioned that we do not meet any type of estimation in
the rate of convergence in any other above mentioned matrix spectral
factorization algorithms.

\subsection{Notation} Let $D:=\{z\in \bC :|z|<1\}$ and $\bT:=\partial D=\{z\in
\bC:|z|=1\}$. $L_p(\bT)$, $p\geq 1$, denotes the class of $p$th
integrable complex functions with usual norm $\|\cdot\|_{L_p}$. For
$f\in L_1(\bT)$, the $n$th Fourier coefficient of $f$ is denoted by
$c_n(f)$.

${\mathcal A}(D)$ denotes the class of analytic functions in $D$.
The Hardy space
$$
H_p=\left\{f\in{\mathcal A}(D): \|f\|_{H_p}=\left(\sup\limits_{r<1}
\int\nolimits_0^{2\pi}|f(re^{it})|^p\,dt<\iy\right)^{\frac
1p}<\iy\right\}
$$
and the space
$$
L_p^+(\bT):= \big\{f\in L_p(\bT): c_n(f)=0\text{ for }n<0\big\}
$$
of the boundary values of  functions from $H_p$ are naturally
identified. The ``$+$" superscript of a function $f^+$ emphasizes
that the function belongs to $L_p^+$. $f\in H_p$ is called outer, we
denote $f\in H_p^O$, if
$$f(z)=c \cdot\exp\left(\frac 1{2\pi}
\int\nolimits_0^{2\pi}\frac{e^{it}+z}{e^{it}-z}\log
|f(e^{it})|\,dt\right),\;\;|c|=1,\;\;z\in D.$$

A matrix function $U(z)$, $z\in \bT$, is called unitary if
$U(z)U^*(z)=I$, where $I$ is the identity matrix of suitable
dimension and $U^*(z)=\ol{U(z)}^T$.

Among the several equivalent norms of a matrix-function we select
the maximum norm of its entries
$$
\|M(z)\|=\max_{i,j}\|M_{ij}(z)\|,
$$
and $M\in L_p$ means that $M_{ij}\in L_p$.

\subsection{The Janashia-Lagvilava algorithm} In order to formulate exactly our result, we
need to describe the Janashia-Lagvilava spectral factorization
method [3] in more details.

Given a positive definite matrix-function
\begin{equation}
S(z)=\begin{pmatrix}a(z)& b(z)\\[2mm] \ol{b(z)}&
    c(z)\end{pmatrix}
\end{equation}
having the properties
\begin{gather}
a(z), b(z), c(z)\in L_1(\bT),\\
a(z)\geq 0,\; a(z)c(z)-|b(z)|^2\geq 0 \text{ for a.a. }z\in\bT,\\
\log \Delta(z):=\log \big(a(z)c(z)-|b(z)|^2\big)\in L_1(\bT),
\end{gather}
the lower-upper triangular factorization is performed at first
\begin{equation}
S(z)=\begin{pmatrix} g^+(z)& 0\\[2mm] \vf(z)& f^+(z)\end{pmatrix}
        \begin{pmatrix} \ol{g^+(z)}& \ol{\vf(z)}\\[2mm] 0
        & \ol{f^+(z)}\end{pmatrix}
\end{equation}
 with the scalar spectral factors
of $a(z)$ and $\Delta(z)/a(z)$, respectively, on the diagonal, i.e.
\begin{gather}
g^+(z)\ol{g^+(z)}=|g^+(z)|^2=a(z),\\f^+(z)\ol{f^+(z)}=|f^+(z)|^2=\frac{\Delta(z)}{a(z)}\,,
\end{gather}
and  $\vf(z)=\ol{b(z)}/\ol{{g^+(z)}}$ in (14). The relations
(11)-(13) imply that $\log a(z)\in L_1(\bT)$ as well, so that the
scalar spectral factorizations in (15) and (16) exist.

The equation $|\vf(z)|^2+|f^+(z)|^2=c(z)$ in (14) and the relations
 (15) and (11) imply that $g^+$, $f^+$ and $\vf$ are square integrable
functions, namely,
\begin{equation}
    \max \{\|g^+\|_{L_2},\|f^+\|_{L_2},\|\vf\|_{L_2}\}\leq
    \sqrt{\|S\|_{L_1}}.
\end{equation}
Thus
\begin{equation}
g^+\,,\,f^+\in H_2^O,\text{ and }\vf\in L_2(\bT).
\end{equation}
For certainty, it is assumed that
\begin{equation}
g^+(0)>0\;\text{ and } f^+(0)>0.
\end{equation}
Next the function $\vf$ is approximated in $L_2$ norm by its Fourier
series
$$
\vf_{n}(z)=\sum_{k=-n}^{\iy} c_k(\vf)z^{n},
$$
so that
\begin{equation}
\|\vf_{n}(z)-\vf(z)\|_{L_2}\to 0,
\end{equation}
and the matrix function $S(z)$ (see (10), (14)) is approximated in
$L_1$ norm by
$$
 S_{n}(z)=\begin{pmatrix} g^+(z)& 0\\[2mm] \vf_{n}(z)& f^+(z)\end{pmatrix}
        \begin{pmatrix} \ol{g^+(z)}& \ol{\vf_{n}(z)}\\[2mm] 0
        & \ol{f^+(z)}\end{pmatrix}.
$$
As it was mentioned in  Introduction, an explicit factorization (6)
is performed in [3] by proving that
\begin{equation}
\chi_{n}^+(z)=\begin{pmatrix} g^+(z)& 0\\[2mm] \vf_{n}(z)& f^+(z)\end{pmatrix}
        \begin{pmatrix}\al_{n}^+(z)&\bt_{n}^+(z)\\[2mm]
        -\ol{\bt_{n}^+(z)}&\ol{\al_{n}^+(z)}\end{pmatrix}
\end{equation}
where the second multiplier in (21) is a unitary matrix-function
with determinant 1,
\begin{equation}
|\al_{n}^+(z)|^2+|\bt_{n}^+(z)|^2=1,\;\; |z|=1,
\end{equation}
and $\al_{n}^+(z)$, $\bt_{n}^+(z)$ are polynomials of order $n$,
$$
\al_{n}^+(z)=\sum_{k=0}^n a_kz^k\;\text{ and
}\bt_{n}^+(z)=\sum_{k=0}^n b_kz^k\,,
$$
whose coefficients can be found by solving a certain system of
linear algebraic equations of order $n$ (see [3], the system (17)
therein).

The equations (21) and (22) provide that (see [3], (21) therein)
\begin{equation}
\det\chi_{n}^+(z)= g^+(z) f^+(z),\;\; z\in D,\;\;n=1,2,\ldots
\end{equation}

Imposing all the same uniqueness conditions on matrices
$\chi_{n}^+(0)$, $n=1,2,\ldots$, say, to be a lower triangular with
positive diagonal entries,
\begin{equation}
[\chi_{n}^+(0)]_{12}=0,\;\;[\chi_{n}^+(0)]_{11}>0,\;\;[\chi_{n}^+(0)]_{22}>0,\;\;
\end{equation}
the convergence (7) is proved (see [3], Sect. 3).

\subsection{The convergence rate estimation}
The closeness of the approximate spectral factor $\chi_{n}^+$ to
$\chi^+$ is estimated in the following
\begin{theorem}
Let $S, S_n,\chi^+,\chi^+_n$ be as in Section 3. Then
\begin{equation}
\|\chi^+_n-\chi^+\|_{L_1}\leq \sqrt{2}\|S\|^{1/2}_{L_1}
\|S^{-1}\|^{1/4}_{L_1}\|\vf_n-\vf\|^{1/2}_{L_2}+
\|\vf_n-\vf\|_{L_1}.
\end{equation}
\end{theorem}
\begin{proof}
We suppress the argument $z$ in what follows. The equation (22)
implies
\begin{equation}
\|\al_n^+\|_{L_\infty}\leq 1 \text{ and }\|\bt_n^+\|_{L_\infty}\leq
1\;\;\;n=1,2,\ldots.
\end{equation}

Taking the product in (21), we get
\begin{equation}
    \chi^+_n= \begin{pmatrix}g^+\al_n^+&g^+\bt_n^+\\[2mm]
    \vf_n\al_n^+-f^+\ol{\bt_n^+}&
    \vf_n\bt_n^++f^+\ol{\al_n^+}\end{pmatrix}=: \begin{pmatrix}g^+\al_n^+&g^+\bt_n^+\\[2mm]
    \Psi_{1n}^+ &\Psi_{2n}^+\end{pmatrix}.
\end{equation}
The claim of the algorithm that $\chi^+_n$ is a spectral factor
implies that
\begin{equation}
\Psi_{1n}^+ \text{ and } \Psi_{2n}^+\in L_2^+.
\end{equation}
The fact that $\Psi_{1n}^+ $ and $\Psi_{2n}^+$ are square integrable
follows also from (18) and (26).

Writing $\chi_m^+$ in the similar form as  $\chi_n^+$ in (27) and
taking the difference, we get
\begin{gather*}
\|\chi_n^+-\chi_m^+\|_{L_1}=\max\big\{\|g^+(\al_n^+-\al_m^+)\|_{L_1},\,\|g^+(\bt_n^+-\bt_m^+)\|_{L_1},\,\\
\|\vf_n(\al_n^+-\al_m^+)+(\vf_n-\vf_m)\al_m^+-f^+(\ol{\bt_n^+}-\ol{\bt_m^+})\|_{L_1},\,\\
\|\vf_n(\bt_n^+-\bt_m^+)+(\vf_n-\vf_m)\bt_n^++f^+(\ol{\al_n^+}-\ol{\al_m^+})\|_{L_1} \big\}
\end{gather*}
\newpage
Thus, by virtue of the Cauchy-Schwartz inequality,
$\|fg\|_{L_1}\leq\|f\|_{L_2}\times\|g\|_{L_2}$, (17), and (26), we
have
\begin{equation}
\|\chi_n^+-\chi_m^+\|_{L_1}\leq \sqrt{\|S\|_{L_1}}
\max\big\{\|\al_n^+-\al_m^+\|_{L_2},\|\bt_n^+-\bt_m^+\|_{L_2}\big\}+
\|\vf_n^+-\vf_m^+\|_{L_1}.
\end{equation}

Consider now the second rows of $\chi^+_n$ and $\chi^+_m$.
\begin{equation}
\begin{cases} \vf_n\al_n^+-f^+\ol{\bt_n^+}
        =\Psi_{1n}^+\,,\\[2mm]
        \vf_n\bt_n^++f^+\ol{\al_n^+}=
        \Psi_{2n}^+\,, \end{cases}\;\;\;
         \begin{cases}
 \vf_m\al_m^+-f^+\ol{\bt_m^+}
        =\Psi_{1m}^+\,,\\[2mm]
        \vf_m\bt_m^++f^+\ol{\al_m^+}=
       \Psi_{2m}^+ \,.\end{cases}
\end{equation}
It follows from (30) that
$$
\begin{cases} \vf_n(\al_n^+-\al_m^+)-f^+(\ol{\bt_n^+}-\ol{\bt_m^+})
        =\Psi_{1n}^+-\Psi_{1m}^+ -(\vf_n-\vf_m)\al_m^+\,,\\[2mm]
        \vf_n(\bt_n^+-\bt_m^+)+f^+(\ol{\al_n^+}-\ol{\al_m^+})=
        \Psi_{2n}^+ -\Psi_{2m}^+-(\vf_n-\vf_m)\bt_m^+\,.\end{cases}
$$
Subtracting the first equation times $(\bt_n^+-\bt_m^+)$ from the
second equation times $(\al_n^+-\al_m^+)$, we get
\begin{gather*}
  f^+\big(|\al_n^+-\al_m^+|^2+|\bt_n^+-\bt_m^+|^2\big)=
(\Psi_{2n}^+ -\Psi_{2m}^+)(\al_n^+-\al_m^+)-\\
(\Psi_{1n}^+-\Psi_{1m}^+)(\bt_n^+-\bt_m^+) -(\vf_n-\vf_m)
(\bt_n^+\al_n^+- \bt_m^+\al_m^+).
\end{gather*}
Hence
\begin{gather}
|\al_n^+-\al_m^+|^2+|\bt_n^+-\bt_m^+|^2= \frac 1{f^+}(\Psi_{2n}^+
-\Psi_{2m}^+)(\al_n^+-\al_m^+)-\\
\frac 1{f^+}(\Psi_{1n}^+-\Psi_{1m}^+)(\bt_n^+-\bt_m^+) - \frac
1{f^+}(\vf_n-\vf_m) (\bt_n^+\al_n^+- \bt_m^+\al_m^+).\notag
\end{gather}
Since
$$
S^{-1}=\frac 1\Delta \begin{pmatrix}c& -b\\[2mm] -\ol{b}& a\end{pmatrix}\in L_1
$$
(see (10), (8)) and $\left|\frac 1{f^+}\right|^2=\frac a\Delta$ (see
(16)), we have
\begin{equation}
\big\|\frac 1{f^+}\big\|_{L_2}^2\leq\|S^{-1}\|_{L_1}<\iy.
\end{equation}
Thus, by virtue of (32), (28), (18), (26), and the Cauchy-Schwartz
inequality, the summands of the right-hand side expression in the
equation (31) are integrable. The generalization of Smirnov's
theorem (see [5], p. 109) claims that if $\Phi=\Phi_1/\Phi_2$, where
$\Phi_1\in H_2$ and $\Phi_2\in H_2^O$, and the boundary values of
$\Phi$ belongs to $L_1(\bT)$, then $\Phi\in H_1$. Thus the first two
functions on the right-hand side of the equation (31) are from
$H_1=L_1^+(\bT)$. Note also that the left-hand side function is
positive. We can write
\begin{gather}
c_0\{|\al_n^+-\al_m^+|^2+|\bt_n^+-\bt_m^+|^2\}=c_0\{\frac 1{f^+}
(\Psi_{2n}^+ -\Psi_{2m}^+)(\al_n^+-\al_m^+)\}-\\c_0\{\frac
1{f^+}(\Psi_{1n}^+-\Psi_{1m}^+)(\bt_n^+-\bt_m^+)\} -c_0\{\frac
1{f^+}(\vf_n-\vf_m) (\bt_n^+\al_n^+- \bt_m^+\al_m^+)\}.\notag
\end{gather}
Now, it follows from (27), (24) and (19) that
\begin{equation}
\al_n^+(0)>0\text{ and }\bt_n^+(0)=0,\;\;\;n=1,2,\ldots,
\end{equation}
so that
\begin{gather}
c_0\{\frac
1{f^+}\big((\Psi_{1n}^+-\Psi_{1m}^+)(\bt_n^+-\bt_m^+)\}=\\\frac
1{f^+(0)}\big((\Psi_{1n}^+(0)-\Psi_{1m}^+(0))(\bt_n^+(0)-\bt_m^+(0))=0.\notag
\end{gather}
Since $
\det(\chi_{{}_n}^+(0))=g^+(0)f^+(0)=g^+(0)\al^+(0)\Psi_{2n}^+(0)$
(see, respectively, (23) and (27), (34)), we have
$\Psi_{2n}^+(0)={f^+(0)}/{\al_n^+(0)}$, and analogously
$\Psi_{2m}^+(0)={f^+(0)}/{\al_m^+(0)}$. Thus, taking into account
(34),
\begin{gather}
c_0\{\frac 1{f^+} (\Psi_{2n}^+ -\Psi_{2m}^+)(\al_n^+-\al_m^+)\}=\\
\frac 1{f^+(0)}(\Psi_{2n}^+
(0)-\Psi_{2m}^+(0))(\al_n^+(0)-\al_m^+(0))=\notag\\
\left(\frac{1}{\al_n^+(0)}
-\frac{1}{\al_m^+(0)}\right)(\al_n^+(0)-\al_m^+(0))=-
\frac{|\al_n^+(0)-\al_m^+(0)|^2}{\al_n^+(0)\al_m^+(0)}\leq 0.\notag
\end{gather}
It follows from (33), (35), (36) and (26) that
\begin{gather*}
c_0\{|\al_n^+-\al_m^+|^2+|\bt_n^+-\bt_m^+|^2\}\leq \\\big|c_0\{\frac
1{f^+}(\vf_n-\vf_m) (\bt_n^+\al_n^+- \bt_m^+\al_m^+)\}\big|\leq
\frac 1\pi \big\|\frac 1{f^+}(\vf_n-\vf_m)\big\|_{L_1}.
\end{gather*}
Hence
\begin{gather*}
\|\al_n^+-\al_m^+\|_{L_2}^2+\|\bt_n^+-\bt_m^+\|_{L_2}^2= \\2\pi
c_0\{|\al_n^+-\al_m^+|^2+|\bt_n^+-\bt_m^+|^2\} \leq 2\|\frac
1{f^+}\|_{L_2}\|(\vf_n-\vf_m)\big\|_{L_2}
\end{gather*}
and, taking
into account (32), we get
\begin{equation}
\max\big\{\|\al_n^+-\al_m^+\|_{L_2},\|\bt_n^+-\bt_m^+\|_{L_2}\big\}\leq
\left(
 2\sqrt{\|S^{-1}\|_{L_1}}\|(\vf_n-\vf_m)\big\|_{L_2}\right)^{\frac12}
\end{equation}
It follows now from (29) and (37) that
$$
\|\chi^+_n-\chi^+_m\|_{L_1}\leq \sqrt{2}\|S\|^{1/2}_{L_1}
\|S^{-1}\|^{1/4}_{L_1}\|\vf_n-\vf_m\|^{1/2}_{L_2}+
\|\vf_n-\vf_m\|_{L_1}
$$
and by virtue of (7) and (20), we get (25).
\end{proof}

\subsection*{Acknowledgements} The authors are grateful to
Professor in Communications Anthony Ephremides (the University of
Maryland) for attracting our attention to the practical importance
of the problem considered in this paper.

We are also obliged to the Japan Society for the Promotion of
Science for the financial support during this research.

\vskip+0.5cm

\noindent Lasha Ephremidze          \hfill Nobuhiko Fujii

\noindent Razmadze Mathematical Institute   \hfill Department of
Mathematics

\noindent Georgian Academy of Sciences     \hfill  Tokai University

\noindent Tbilisi 0193, Georgia     \hfill   Shizuoka  424-8610,
Japan

\noindent E-mail: lasha@rmi.acnet.ge    \hfill  E-mail:
nfujii@scc.u-tokai.ac.jp

\end{document}